\documentclass[12pt,reqno,a4paper]{amsart}
\usepackage{mathrsfs}
\usepackage{latexsym,amsmath,amsfonts,amssymb,amsthm}
\theoremstyle{plain}
\newtheorem{Thm}{Theorem}
\newtheorem{Lem}{Lemma}

\newtheorem{Conj}{Conjecture}
\theoremstyle{definition}

\theoremstyle{remark}

\def\Z{\mathbb Z}

\def\D{\mathscr{D}}
\def\N{\mathbb N}

\def\1{{\bf 1}}
\def\floor #1{\left\lfloor{#1}\right\rfloor}

\def\qbinom #1#2#3{\genfrac{[}{]}{0pt}{}{#1}{#2}_{#3}}
\def\pmod #1{\ ({{\rm mod}}\ #1)}

\begin{document}

\title{Factors of alternative binomials sums}
\author{Hui-Qin Cao}
\address{Department of Mathematics, Nanjing Audit University,
Nanjing 210029, People's Republic of China} \email{caohq@nau.edu.cn}
\author{Hao Pan}
\address{Department of Mathematics, Nanjing University,
Nanjing 210093, People's Republic of China}
\email{haopan79@yahoo.com.cn}
\thanks{The first author was supported by
the National Natural Science Foundation of China (Grant No.
10871087). And the second author was supported by
the National Natural Science Foundation of China (Grant No.
10771135).}\keywords{}
\date{}
\maketitle

\begin{abstract}
We confirm several conjectures of Guo, Jouhet and Zeng concerning
the factors of alternative binomials sums.
\end{abstract}

\section{Introduction}
\setcounter{equation}{0} \setcounter{Thm}{0} \setcounter{Lem}{0}
\setcounter{Cor}{0}

It is well-known that
$$
\sum_{k=0}^{n}(-1)^k\binom{n}{k}=(1-1)^{n}=0
$$
for every positive integer $n$. However, there are two unfamiliar identities in the same flavor \cite[Eqs. (3.81) and (6.6)]{Gould72}:
\begin{equation}
\label{be1}
\sum_{k=0}^{2n}(-1)^k\binom{2n}{k}^2=(-1)^n\binom{2n}{n}
\end{equation}
and
\begin{equation}
\label{be2}
\sum_{k=0}^{2n}(-1)^k\binom{2n}{k}^3=(-1)^n\binom{2n}{n}\binom{3n}{n}
\end{equation}
for any $n\geq 1$. Unfortunately, by using asymptotic methods, de Bruijn \cite{Bruijn81} has showed that no closed
form exists for the sum $\sum_{k=0}^n(-1)^k\binom{n}{k}^a$ when $a\geq4$. Observe that the right sides of (\ref{be1}) and (\ref{be2})
are both divisible by $\binom{2n}{n}$.
Motivated by (\ref{be1}) and (\ref{be2}), in \cite{Calkin98}, Calkin established the following interesting congruence:
\begin{equation}
\label{Calkin}
\sum_{k=0}^{2n}(-1)^k\binom{2n}{k}^r\equiv 0\pmod{\binom{2n}{n}}
\end{equation} for any positive integers $n$ and $r$. Nine years later,
Guo, Jouhet and Zeng \cite{GuoJouhetZeng07} generalized Calkin's result and showed that for
any positive integers $n_1,\ldots,n_h,n_{h+1}=n_1$,
\begin{equation}
\label{gjz}
\sum_{k=-n_1}^{n_1}(-1)^k\prod_{i=1}^h\binom{n_i+n_{i+1}}{n_i+k}\equiv
0\pmod{\binom{n_1+n_h}{n_1}}
\end{equation}
In fact, they proved a $q$-analogue of (\ref{gjz}):
\begin{equation}
\label{gjzq}
\sum_{k=-n_1}^{n_1}(-1)^kq^{\binom{k}{2}}\prod_{i=1}^h\qbinom{n_i+n_{i+1}}{n_i+k}q\equiv
0\pmod{\qbinom{n_1+n_r}{n_1}q},
\end{equation}
where the above congruence is considered over the polynomials ring $\Z[q]$.

Based on some computer experiments, Guo, Jouhet and Zeng proposed  several conjectures on alternative binomial sums:
\begin{Conj}
\label{Conj1} For any positive integers $m$ and $n$,
\begin{equation}
\gcd\bigg(\sum_{k=0}^{2n}(-1)^k\binom{2n}{k}^r:\,
r=m,m+1,\ldots\bigg)=\binom{2n}n,
\end{equation}
where $\gcd(a_1,a_2,\ldots)$ denotes the greatest common divisor
of $a_1,a_2,\ldots$.
\end{Conj}
\begin{Conj}
\label{Conj2}
For any positive integers $r,s,t$ and $n$,
\begin{equation}
\label{cj2c1}
\sum_{k=-n}^{n}(-1)^k\binom{6n}{3n+k}^r\binom{4n}{2n+k}^s\binom{2n}{n+k}^t\equiv
0\pmod{2\binom{6n}{n}},
\end{equation}
\begin{equation}
\label{cj2c2}
\sum_{k=-n}^{n}(-1)^k\binom{6n}{3n+k}^r\binom{4n}{2n+k}^s\binom{2n}{n+k}^t\equiv
0\pmod{6\binom{6n}{3n}}.
\end{equation}
Furthermore, if $(r,s,t)\not=(1,1,1)$, then
\begin{equation}
\label{cj2c3}
\sum_{k=-n}^{n}(-1)^k\binom{8n}{4n+k}^r\binom{4n}{2n+k}^s\binom{2n}{n+k}^t\equiv
0\pmod{2\binom{8n}{3n}}.
\end{equation}
\end{Conj}

In this paper, we shall confirm these conjectures. For a prime $p$
and an integer $n$, let $\nu_p(n)$ denote the greatest integer such
that $p^{\nu_p(n)}\mid n$. In particular, we set $\nu_p(0)=+\infty$.
Let $\phi$ denote the Euler totient function. Clearly Conjecture
\ref{Conj1} is implied by the following theorem.
\begin{Thm}
\label{thm1}
Suppose that $n$ is a positive integer and $r$ is a positive integer with $r\equiv 2\pmod{\phi({\binom{2n}n}){\binom{2n}n}}$. Then
$$
\nu_p\bigg(\sum_{k=0}^{2n}(-1)^k\binom{2n}{k}^r\bigg)=\nu_p\bigg(\binom{2n}n\bigg)
$$
for each prime divisor $p$ of $\binom{2n}n$.
\end{Thm}

For a positive integer $n$, define
$$
[n]_q=\frac{1-q^n}{1-q}=1+q+q^2+\cdots+q^{n-1}.
$$
And define the $q$-binomial coefficient
$$
\qbinom{n}{k}q=\begin{cases}\prod_{j=1}^k\frac{1-q^{n+1-j}}{1-q^j},\qquad&\text{if }n\geq k\geq 1\\
1,\qquad&\text{if }k=0,\\
0,\qquad&\text{if }k<0\text{ or }n<k.
\end{cases}
$$
Applying (\ref{gjzq}), it is not difficult (see \cite[Theorem 4.7,
Corollary 4.10 and Corollary 4.11]{GuoJouhetZeng07}) to deduce that
\begin{equation}
\label{cj2c1q}
\sum_{k=-n}^{n}(-1)^kq^{\binom{k}2}\qbinom{6n}{3n+k}q^r\qbinom{4n}{2n+k}q^s\qbinom{2n}{n+k}q^t\equiv
0\pmod{\qbinom{6n}{n}q},
\end{equation}
\begin{equation}
\label{cj2c2q}
\sum_{k=-n}^{n}(-1)^kq^{\binom{k}2}\qbinom{6n}{3n+k}q^r\qbinom{4n}{2n+k}q^s\qbinom{2n}{n+k}q^t\equiv
0\pmod{\qbinom{6n}{3n}q},
\end{equation}
and
\begin{equation}
\label{cj2c3q}
\sum_{k=-n}^{n}(-1)^kq^{\binom{k}2}\qbinom{8n}{4n+k}q^r\qbinom{4n}{2n+k}q^s\qbinom{2n}{n+k}q^t\equiv
0\pmod{\qbinom{8n}{3n}q}.
\end{equation}
Now we shall prove that
\begin{Thm}
\label{thm2} Let $\alpha=\nu_2(n)$ and $\beta=\nu_3(n)$. For positive integers $r,s,t$,
\begin{equation}
\label{t2c1}
{\qbinom{6n}{n}q}^{-1}\sum_{k=-n}^{n}(-1)^kq^{\binom{k}2}\qbinom{6n}{3n+k}q^r\qbinom{4n}{2n+k}q^s\qbinom{2n}{n+k}q^t\equiv
0\pmod{[2]_{q^{2^\alpha}}}
\end{equation}
and
\begin{equation}
\label{t2c2}
{\qbinom{6n}{3n}q}^{-1}\sum_{k=-n}^{n}(-1)^kq^{\binom{k}2}\qbinom{6n}{3n+k}q^r\qbinom{4n}{2n+k}q^s\qbinom{2n}{n+k}q^t\equiv
0\pmod{[2]_{q^{2^\alpha}}[3]_{q^{2^\alpha}}}.
\end{equation}
Further, we have
\begin{align}
\label{t2c3}
&{\qbinom{8n}{3n}q}^{-1}\sum_{k=-n}^{n}(-1)^kq^{\binom{k}2}\qbinom{8n}{4n+k}q^r\qbinom{4n}{2n+k}q^s\qbinom{2n}{n+k}q^t\notag\\
\equiv&
\begin{cases}0\pmod{[2]_{q^{2^\alpha}}},&\qquad\text{if }t\geq 2,\\
0\pmod{[2]_{q^{2^{\alpha+1}}}},&\qquad\text{if }s\geq 2,\text{ or }r\geq 2\text{ and }n\equiv 3\cdot2^\alpha\pmod{2^{\alpha+2}},\\
0\pmod{[2]_{q^{2^{\alpha+2}}}},&\qquad\text{if }r\geq 2\text{ and
}n\equiv 2^\alpha\pmod{2^{\alpha+2}}.
\end{cases}
\end{align}
\end{Thm}
Let us explain why Theorem \ref{thm2} implies Conjecture\ref{Conj2}.
For example, since $[2]_{q^{2^\alpha}}$ is a primitive polynomial (a
polynomial with integral coefficients is called primitive if the
greatest common divisor of its coefficients is 1), by (\ref{t2c1}),
there exists a polynomial $H(q)$ with integral coefficients such
that
$$
\sum_{k=-n}^{n}(-1)^kq^{\binom{k}2}\qbinom{6n}{3n+k}q^r\qbinom{4n}{2n+k}q^s\qbinom{2n}{n+k}q^t=H(q)[2]_{q^{2^\alpha}}{\qbinom{6n}{n}q}.
$$
Thus substituting $q=1$ in the above equation, we get
$$
\sum_{k=-n}^{n}(-1)^k\binom{6n}{3n+k}^r\binom{4n}{2n+k}^s\binom{2n}{n+k}^t=2H(1){\binom{6n}{n}},
$$
that is,
$$
\sum_{k=-n}^{n}(-1)^k\binom{6n}{3n+k}^r\binom{4n}{2n+k}^s\binom{2n}{n+k}^t\equiv0\pmod{2{\binom{6n}{n}}}.
$$

The proofs of Theorems \ref{thm1} and \ref{thm2} will be proposed in Sections 2 and 3.

\section{Proof of Theorem \ref{thm1}}
\setcounter{equation}{0} \setcounter{Thm}{0} \setcounter{Lem}{0}
\setcounter{Cor}{0}

Suppose that $p$ is an arbitrary prime divisor of $\binom{2n}n$ and
$\nu_p(\binom{2n}n)=\gamma$. Suppose that $r>2$ be an integer such
that
$$
r\equiv 2\pmod{\phi(p^{\gamma+1})}.
$$
It is easy to see that $r\geq\gamma+1$. Then
$$
\sum_{k=0}^{2n}(-1)^k\binom{2n}{k}^r\equiv\sum_{\substack{0\leq
k\leq 2n\\
p\nmid\binom{2n}k}}(-1)^k\binom{2n}{k}^2\pmod{p^{\gamma+1}}.
$$
Thus Theorem \ref{thm1} easily follows from:
\begin{Lem} Let $p$ be a prime and $n$ be a positive integer. Then
\begin{equation}
\nu_p\bigg(\sum_{\substack{0\leq k\leq 2n\\ p\nmid\binom{2n}k}}(-1)^k\binom{2n}{k}^2\bigg)=\nu_p\bigg(\binom{2n}{n}\bigg).
\end{equation}
\end{Lem}

Notice that
$$
\sum_{\substack{0\leq k\leq 2n\\ p\nmid\binom{2n}k}}(-1)^k\binom{2n}{k}^2+\sum_{\substack{0\leq k\leq 2n\\ p\mid\binom{2n}k}}(-1)^k\binom{2n}{k}^2=\sum_{k=0}^{2n}(-1)^k\binom{2n}{k}^2=(-1)^n\binom{2n}n.
$$
So we only need to prove that
\begin{Lem} For each $r\geq 1$,
\begin{equation}
\nu_p\bigg(\sum_{\substack{0\leq k\leq 2n\\ p\mid\binom{2n}k}}(-1)^k\binom{2n}{k}^r\bigg)\geq r-1+\nu_p\bigg(\binom{2n}{n}\bigg).
\end{equation}
\end{Lem}

Let
$$
\D_{n,k}=\{d\in\N:\, \floor{n/d}>\floor{k/d}+\floor{(n-k)/d}\},
$$
where $\floor{x}=\max\{z\in\Z:\,z\leq x\}$. Note that $p\mid\binom{2n}k$ if and only if the set $\{\beta:\, p^\beta\in\D_{2n,k}\}$ is non-empty.
Letting $h=\floor{\log_p(2n)}+1$, we have
\begin{align*}
\sum_{\substack{0\leq k\leq 2n\\ p\mid\binom{2n}k}}(-1)^k\binom{2n}{k}^2=&
\sum_{k=0}^{2n}(-1)^k\binom{2n}{k}^2\sum_{\emptyset\not=I\subseteq\{\alpha:\, p^\alpha\in\D_{2n,k}\}}(-1)^{|I|-1}\\
=&\sum_{\emptyset\not=I\subseteq\{1,2,\ldots,h\}}(-1)^{|I|-1}\sum_{\substack{0\leq k\leq 2n\\ p^\alpha\in\D_{2n,k},\ \forall\alpha\in I}}(-1)^k\binom{2n}{k}^2.
\end{align*}
Hence it suffices to show that
\begin{Lem} For each $\emptyset\not=I\subseteq\{1,\ldots,h\}$,
\begin{equation}
\label{gcdlc3}
\nu_p\bigg(\sum_{\substack{0\leq k\leq 2n\\ p^\alpha\in\D_{2n,k},\ \forall\alpha\in I}}(-1)^k\binom{2n}{k}^r\bigg)\geq(r-1)|I|+\nu_p\bigg(\binom{2n}{n}\bigg).
\end{equation}
\end{Lem}

It is not difficult to see that
$$
\qbinom{n}{k}q=\prod_{d\in\D_{n,k}}\Phi_d(q),
$$
where $\Phi_d(q)$ is the $d$-th cyclotomic polynomial. In
particular, we have
$$
\Phi_{p^\alpha}(q)=\frac{1-q^{p^\alpha}}{1-q^{p^{\alpha-1}}}=[p]_{q^{p^{\alpha-1}}}
$$
for every prime $p$ and integer $\alpha\geq 1$. Thus (\ref{gcdlc3})
is an immediate consequence of the following $q$-congruence.
\begin{Lem}
\begin{equation}
\label{gcdlc4}
\sum_{\substack{0\leq k\leq 2n\\ p^\alpha\in\D_{2n,k},\ \forall\alpha\in I}}(-1)^kq^{\binom{k}2}\qbinom{2n}{k}q^r\equiv 0\pmod{\prod_{\alpha\in I}\Phi_{p^{\alpha}}(q)^r\prod_{\substack{\beta\not\in I\\ p^\beta\in\D_{2n,n}}}\Phi_{p^{\beta}}(q)}.
\end{equation}
\end{Lem}
\begin{proof}
We need a $q$-analogue of well-known Lucas' congruence (cf. \cite{Sagan92}):
\begin{equation}
\label{qLucas}
\qbinom{x_1d+x_2}{y_1d+y_2}q\equiv\binom{x_1}{y_1}\qbinom{x_2}{y_2}q\pmod{\Phi_d(q)}
\end{equation}
for every $d\geq 2$, where $0\leq x_2,y_2<d$.

For any $\beta$ with $\beta\not\in I$ and $p^\beta\in\D_{2n,n}$, write
$n=n_1p^\beta+n_2$ with $0\leq n_2<p^\beta$. Since $p^\beta\in\D_{2n,n}$, we have $2n_2\geq p^\beta$. For any $k=k_1p^\beta+k_2$ with $0\leq k_2<p^\beta$, by (\ref{qLucas}),
$$
\qbinom{2n}{k}q\equiv\binom{2n_1+1}{k_1}\qbinom{2n_2-p^\beta}{k_2}q\pmod{\Phi_{p^\beta}(q)}.
$$
Hence
$$
\qbinom{2n}{k}q\equiv0\pmod{\Phi_{p^\beta}(q)}.
$$
provided that $2n_2-p^\beta<k_2$.

Suppose that $2n_2-p^\beta\geq k_2$. Assume that
$I=\{\alpha_1,\alpha_2,\ldots,\alpha_u\}$ with
$$
\alpha_1<\alpha_2<\ldots<\alpha_v<\beta<\alpha_{v+1}<\ldots<\alpha_{u}.
$$
When $1\leq j\leq v$, we have
\begin{align*}
&\floor{\frac{2n}{p^{\alpha_j}}}-\floor{\frac{k}{p^{\alpha_j}}}-\floor{\frac{2n-k}{p^{\alpha_j}}}\\
=&
\floor{\frac{(2n_1+1)p^\beta+2n_2-p^\beta}{p^{\alpha_j}}}-\floor{\frac{k_1p^\beta+k_2}{p^{\alpha_j}}}-\floor{\frac{(2n_1+1-k_1)p^\beta+2n_2-p^\beta-k_2}{p^{\alpha_j}}}\\
=&\floor{\frac{2n_2-p^\beta}{p^{\alpha_j}}}-\floor{\frac{k_2}{p^{\alpha_j}}}-\floor{\frac{2n_2-p^\beta-k_2}{p^{\alpha_j}}}.
\end{align*}
It follows that
$p^{\alpha_j}\in\D_{2n,k}$ if and only if $p^{\alpha_j}\in\D_{2n_2-p^\beta,k_2}$ for $1\leq j\leq v$. Similarly,
\begin{align*}
&\floor{\frac{(2n_1+1)p^\beta+2n_2-p^\beta}{p^{\alpha_j}}}-\floor{\frac{k_1p^\beta+k_2}{p^{\alpha_j}}}-\floor{\frac{(2n_1+1-k_1)p^\beta+2n_2-p^\beta-k_2}{p^{\alpha_j}}}\\
=&\floor{\frac{2n_1+1}{p^{\alpha_j-\beta}}}-\floor{\frac{k_1}{p^{\alpha_j-\beta}}}-\floor{\frac{2n_1+1-k}{p^{\alpha_j-\beta}}}
\end{align*}
provided that $\alpha_j>\beta$. Therefore $p^{\alpha_j}\in\D_{2n,k}$
if and only if $p^{\alpha_j-\beta}\in\D_{2n_1+1,k_1}$ for $v+1\leq
j\leq u$. Thus
\begin{align*}
&\sum_{\substack{0\leq k\leq 2n\\ p^\alpha\in\D_{2n,k},\ \forall\alpha\in I}}(-1)^kq^{\binom k2}\qbinom{2n}{k}q^r\\
\equiv& \sum_{\substack{0\leq k_1\leq 2n_1+1\\
p^{\alpha_j-\beta}\in\D_{2n_1+1,k_1},\\ \forall
j\in\{v+1,\ldots,u\}}}(-1)^{k_1p^\beta}q^{\binom{k_1p^\beta}2}\binom{2n_1+1}{k_1}^r
\cdot\sum_{\substack{0\leq k_2\leq 2n_2-p^\beta\\ p^{\alpha_j}\in\D_{2n_2-p^\beta,k_2},\\ \forall j\in\{1,\ldots,v\}}}(-1)^{k_2}q^{\binom{k_2}2}\qbinom{2n_2-p^\beta}{k_2}q^r\\
&\pmod{\Phi_{p^\beta}(q)},
\end{align*}
by noting that
$$
q^{\binom
k2}=q^{\binom{k_1p^\beta+k_2}2}=q^{\binom{k_1p^\beta}2+\binom{k_2}2+k_1k_2p^\beta}\equiv
q^{\binom{k_1p^\beta}2+\binom{k_2}2}\pmod{\Phi_{p^\beta}(q)}.
$$

If $p$ is an odd prime, then
$$
q^{\binom{k_1p^\beta}2}=(q^{p^\beta})^{\frac{k_1(k_1p^\beta-1)}{2}}\equiv1\pmod{\Phi_{p^\beta}(q)}.
$$
And if $p=2$, then we have
$$
q^{\binom{k_12^\beta}2}=(q^{2^{\beta-1}})^{k_1(k_12^\beta-1)}\equiv(-1)^{k_1}\pmod{\Phi_{2^\beta}(q)}
$$
since $1+q^{2^{\beta-1}}=[2]_{q^{2^{\beta-1}}}=\Phi_{2^\beta}(q)$.
Notice that $\D_{2n_1+1,k_1}=\D_{2n_1+1,2n_1+1-k_1}$. We have
\begin{align*}
&\sum_{\substack{0\leq k_1\leq 2n_1+1\\ p^{\alpha_j-\beta}\in\D_{2n_1+1,k_1},\\ \forall j\in\{v+1,\ldots,u\}}}(-1)^{k_1p^\beta}q^{\binom{k_1p^\beta}2}\binom{2n_1+1}{k_1}^r\\
\equiv& \frac{1}{2}\sum_{\substack{0\leq k_1\leq 2n_1+1\\
p^{\alpha_j-\beta}\in\D_{2n_1+1,k_1},\\ \forall
j\in\{v+1,\ldots,u\}}}\big((-1)^{k_1}+(-1)^{2n_1+1-k_1}\big)\binom{2n_1+1}{k_1}^r=0\pmod{\Phi_{p^\beta}(q)}.
\end{align*}

Finally, clearly
\begin{align*}
\sum_{\substack{0\leq k\leq 2n\\ p^\alpha\in\D_{2n,k},\ \forall\alpha\in I}}(-1)^kq^{\binom k2}\qbinom{2n}{k}q^r\equiv0\pmod{\Phi_{p^\alpha}(q)^r}
\end{align*}
for any $\alpha\in I$.
\end{proof}

\section{Proof of Theorem \ref{thm2}}
\setcounter{equation}{0} \setcounter{Thm}{0} \setcounter{Lem}{0}
\setcounter{Cor}{0}

Recalling that $\qbinom{n}{k}q=\prod_{d\in\D_{n,k}}\Phi_d(q)$ and
$\Phi_{p^\alpha}(q)=[p]_{q^{p^{\alpha-1}}}$. Let $\alpha=\nu_2(n)$.
For any $k$ with $\nu_2(k)\not=\alpha$, since
$$
2n\equiv 0\pmod{2^{\alpha+1}}\qquad\text{and}\qquad n+k\not\equiv 0\pmod{2^{\alpha+1}},
$$
we have
$$
\qbinom{2n}{n+k}q\equiv 0\pmod{\Phi_{2^{\alpha+1}}(q)}.
$$
Similarly,
$$
\qbinom{6n}{3n+k}q\equiv 0\pmod{\Phi_{2^{\alpha+1}}(q)}.
$$
Hence
\begin{equation}
\label{c6nn1}
\sum_{\substack{-n\leq k\leq n\\ \nu_2(k)\not=\alpha}}(-1)^kq^{\binom
k2}\qbinom{6n}{3n+k}q^r\qbinom{4n}{2n+k}q^s\qbinom{2n}{n+k}q^t\equiv0\pmod{\Phi_{2^{\alpha+1}}(q)^2}.
\end{equation}
On the other hand, obviously
\begin{align*}
&\sum_{\substack{-n\leq k\leq n\\ \nu_2(k)=\alpha}}(-1)^kq^{\binom
k2}\qbinom{6n}{3n+k}q^r\qbinom{4n}{2n+k}q^s\qbinom{2n}{n+k}q^t\\=&\sum_{\substack{k>0\\ \nu_2(k)=\alpha}}(-1)^kq^{\binom
k2}(1+q^k)\qbinom{6n}{3n+k}q^r\qbinom{4n}{2n+k}q^s\qbinom{2n}{n+k}q^t.
\end{align*}
For any $k$ with $\nu_2(k)=\alpha$, we have
$$
4n\equiv 0\pmod{2^{\alpha+1}}\qquad\text{and}\qquad 2n+k\equiv 2^\alpha\pmod{2^{\alpha+1}},
$$
whence
$$
\qbinom{4n}{2n+k}q\equiv 0\pmod{\Phi_{2^{\alpha+1}}(q)}.
$$
And $1+q^k$ is divisible by $1+q^{2^\alpha}=\Phi_{2^{\alpha+1}}(q)$,
since $k/2^\alpha$ is odd. Thus
\begin{equation}
\label{c6nn2}
\sum_{\substack{-n\leq k\leq n\\ \nu_2(k)=\alpha}}(-1)^kq^{\binom
k2}\qbinom{6n}{3n+k}q^r\qbinom{4n}{2n+k}q^s\qbinom{2n}{n+k}q^t\equiv0\pmod{\Phi_{2^{\alpha+1}}(q)^2}.
\end{equation}
Combining (\ref{c6nn1}) and (\ref{c6nn2}), we have
\begin{equation}
\label{c6nn3}
\sum_{k=-n}^n(-1)^kq^{\binom
k2}\qbinom{6n}{3n+k}q^r\qbinom{4n}{2n+k}q^s\qbinom{2n}{n+k}q^t\equiv0\pmod{\Phi_{2^{\alpha+1}}(q)^2}.
\end{equation}
And by (\ref{c6nn3}) and (\ref{cj2c1q}), we conclude that
\begin{equation*}
\label{c6nn4}
\sum_{k=-n}^n(-1)^kq^{\binom
k2}\qbinom{6n}{3n+k}q^r\qbinom{4n}{2n+k}q^s\qbinom{2n}{n+k}q^t\equiv0\pmod{\Phi_{2^{\alpha+1}}(q)\qbinom{6n}{n}q},
\end{equation*}
since $\Phi_{2^{\alpha+1}}(q)^2\nmid\qbinom{6n}{n}q$.

Let $\beta=\nu_3(n)$. If $\nu_3(k)\leq\beta$, then
$$
6n\equiv 3n\equiv 0\pmod{3^{\beta+1}}\qquad\text{and}\qquad 3n+k\not\equiv 0\pmod{3^{\beta+1}},
$$
whence
$$
\qbinom{6n}{3n+k}q\equiv 0\pmod{\Phi_{3^{\beta+1}}(q)}.
$$
Suppose that $\nu_3(k)>\beta$. If $n\equiv 3^\beta\pmod{3^{\beta+1}}$. Then
$$
4n\equiv3^\beta\pmod{3^{\beta+1}}\qquad\text{and}\qquad 2n+k\equiv 2\cdot3^\beta\pmod{3^{\beta+1}}.
$$
Thus
$$
\qbinom{4n}{2n+k}q\equiv 0\pmod{\Phi_{3^{\beta+1}}(q)}.
$$
And if $n\equiv 2\cdot3^\beta\pmod{3^{\beta+1}}$, then
$$
2n\equiv3^\beta\pmod{3^{\beta+1}}\qquad\text{and}\qquad n+k\equiv 2\cdot3^\beta\pmod{3^{\beta+1}},
$$
whence
$$
\qbinom{2n}{n+k}q\equiv 0\pmod{\Phi_{3^{\beta+1}}(q)}.
$$
This concludes that
\begin{equation}
\label{c6n3n1}
\sum_{k=-n}^n(-1)^kq^{\binom
k2}\qbinom{6n}{3n+k}q^r\qbinom{4n}{2n+k}q^s\qbinom{2n}{n+k}q^t\equiv 0\pmod{\Phi_{3^{\beta+1}}(q)}.
\end{equation}
Since $6n\equiv 3n\equiv 0\pmod{3^{\beta+1}}$,
$3^{\beta+1}\not\in\D_{6n,3n}$, i.e.,
$\Phi_{3^{\beta+1}}(q)\nmid\qbinom{6n}{3n}q$. Thus combining
(\ref{c6nn3}), (\ref{c6n3n1}) and (\ref{cj2c2q}), we get
(\ref{t2c2}).

Finally, let us turn to (\ref{cj2c3}). Suppose that $\nu_2(n)=\alpha$. Since $(r,s,t)\not=(1,1,1)$, we may consider the following three cases:

\medskip\noindent{\bf Case 1}: $t\geq 2$. If $\nu_2(k)\not=\alpha$, then
$$
2n\equiv0\pmod{2^{\alpha+1}}\qquad\text{and}\qquad n+k\not\equiv 0\pmod{2^{\alpha+1}},
$$
whence
$$
\qbinom{2n}{n+k}q\equiv 0\pmod{\Phi_{2^{\alpha+1}}(q)}.
$$
And if $\nu_2(k)=\alpha$, then
$$
8n\equiv4n\equiv0\pmod{2^{\alpha+1}}\qquad\text{and}\qquad 4n+k\equiv 2n+k\equiv 2^\alpha\pmod{2^{\alpha+1}}.
$$
So
$$
\qbinom{8n}{4n+k}q\equiv\qbinom{4n}{2n+k}q\equiv 0\pmod{\Phi_{2^{\alpha+1}}(q)}.
$$
Hence
\begin{equation}
\label{c8n3nc1}
\sum_{k=-n}^n(-1)^kq^{\binom
k2}\qbinom{8n}{4n+k}q^r\qbinom{4n}{2n+k}q^s\qbinom{2n}{n+k}q^t\equiv 0\pmod{\Phi_{2^{\alpha+1}}(q)^2}.
\end{equation}

\medskip\noindent{\bf Case 2}: $s\geq 2$. If $\nu_2(k)\not=\alpha+1$, then
$$
4n\equiv0\pmod{2^{\alpha+2}}\qquad\text{and}\qquad 2n+k\not\equiv 0\pmod{2^{\alpha+2}},
$$
whence
$$
\qbinom{4n}{2n+k}q\equiv 0\pmod{\Phi_{2^{\alpha+2}}(q)}.
$$
Assume that $\nu_2(k)=\alpha+1$. Then
$$
8n\equiv0\pmod{2^{\alpha+2}}\qquad\text{and}\qquad 4n+k\equiv 2^{\alpha+1}\pmod{2^{\alpha+2}}.
$$
It follows that
$$
\qbinom{8n}{4n+k}q\equiv0\pmod{\Phi_{2^{\alpha+2}}(q)}.
$$
And $\Phi_{2^{\alpha+2}}(q)=1+q^{2^{\alpha+1}}$ divides $1+q^k$ since $k/2^{\alpha+1}$ is odd.
Thus
\begin{align}
\label{c8n3nc2}
&\sum_{k=-n}^n(-1)^kq^{\binom
k2}\qbinom{8n}{4n+k}q^r\qbinom{4n}{2n+k}q^s\qbinom{2n}{n+k}q^t\notag\\
\equiv&\sum_{\substack{-n\leq k\leq n\\ \nu_2(k)=\alpha+1}}(-1)^kq^{\binom
k2}\qbinom{8n}{4n+k}q^r\qbinom{4n}{2n+k}q^s\qbinom{2n}{n+k}q^t\notag\\
=&\sum_{\substack{0<k\leq n\\ \nu_2(k)=\alpha+1}}(-1)^kq^{\binom
k2}(1+q^k)\qbinom{8n}{4n+k}q^r\qbinom{4n}{2n+k}q^s\qbinom{2n}{n+k}q^t\notag\\
\equiv&0\pmod{\Phi_{2^{\alpha+2}}(q)^2}.
\end{align}

\medskip\noindent{\bf Case 3}: $r\geq 2$. We consider two subcases:

\medskip\noindent(i) $n\equiv 2^\alpha\pmod{2^{\alpha+2}}$. For any $k$ with $\nu_2(k)\not=\alpha+2$, we have
$$
8n\equiv0\pmod{2^{\alpha+3}}\qquad\text{and}\qquad 4n+k\not\equiv 0\pmod{2^{\alpha+3}}.
$$
So
$$
\qbinom{8n}{4n+k}q\equiv 0\pmod{\Phi_{2^{\alpha+3}}(q)}.
$$
And for any $k$ with $\nu_2(k)=\alpha+2$, we have
$$
4n\equiv2^{\alpha+2}\pmod{2^{\alpha+3}}\qquad\text{and}\qquad 2n+k\equiv2^{\alpha+2}+2^{\alpha+1}\pmod{2^{\alpha+3}}.
$$
Then
$$
\qbinom{4n}{2n+k}q\equiv1+q^k\equiv0\pmod{\Phi_{2^{\alpha+3}}(q)}.
$$
Thus
\begin{align}
\label{c8n3nc3}
&\sum_{k=-n}^n(-1)^kq^{\binom
k2}\qbinom{8n}{4n+k}q^r\qbinom{4n}{2n+k}q^s\qbinom{2n}{n+k}q^t\notag\\
\equiv&\sum_{\substack{0<k\leq n\\ \nu_2(k)=\alpha+2}}(-1)^kq^{\binom
k2}(1+q^k)\qbinom{8n}{4n+k}q^r\qbinom{4n}{2n+k}q^s\qbinom{2n}{n+k}q^t\notag\\
\equiv&0\pmod{\Phi_{2^{\alpha+3}}(q)^2}.
\end{align}

\medskip\noindent(ii) $n\equiv 3\cdot 2^\alpha\pmod{2^{\alpha+2}}$. For any $k$ with $\nu_2(k)<\alpha+2$, we have
$$
8n\equiv4n\equiv0\pmod{2^{\alpha+2}}\qquad\text{and}\qquad 4n+k\not\equiv 0\pmod{2^{\alpha+2}},
$$
whence
$$
\qbinom{8n}{4n+k}q\equiv 0\pmod{\Phi_{2^{\alpha+2}}(q)}.
$$
If $\nu_2(k)\geq\alpha+2$, then
\begin{align*}
4n\equiv0\pmod{2^{\alpha+2}},\ 2n\equiv 2n+k\equiv2^{\alpha+1}\pmod{2^{\alpha+2}}
,\ n+k\equiv3\cdot2^{\alpha}\pmod{2^{\alpha+2}}.
\end{align*}
Hence
$$
\qbinom{4n}{2n+k}q\equiv\qbinom{2n}{n+k}q\equiv0\pmod{\Phi_{2^{\alpha+2}}(q)}
$$
and
\begin{equation}
\label{c8n3nc4}
\sum_{k=-n}^n(-1)^kq^{\binom
k2}\qbinom{8n}{4n+k}q^r\qbinom{4n}{2n+k}q^s\qbinom{2n}{n+k}q^t
\equiv0\pmod{\Phi_{2^{\alpha+2}}(q)^2}.
\end{equation}
From (\ref{c8n3nc1})-(\ref{c8n3nc4}) and (\ref{cj2c3q}),
(\ref{t2c3}) is concluded.

\end{document}